\documentclass[12pt]{article}

\usepackage{amsfonts, amsmath, amsthm, amssymb, tikz, url}
\usepackage[in]{fullpage}
\usepackage[english]{babel}
\usepackage[utf8]{inputenc}
\usepackage{authblk}
\usepackage{float}
\usepackage{subfigure}
\usepackage{graphicx}

\usetikzlibrary{shapes.geometric}

\title{The Diagonal of a Regular Pentagon}
\author{Dan McQuillan}
\affil{Department of Mathematics, Norwich University, Northfield VT 05663\\
\tt dmcquill@norwich.edu}
\author{Dan McQuillan\\               
}

\begin{document}
\maketitle
\abstract{
We demonstrate the well-known fact that the diagonal of a unit regular pentagon has length the golden ratio. While similar to other proofs, we provide a novel figure that we believe is particularly well-suited to achieve full understanding at a glance.}

\begin{figure}[htbp]
   \centering
    \includegraphics[width=0.7\textwidth]{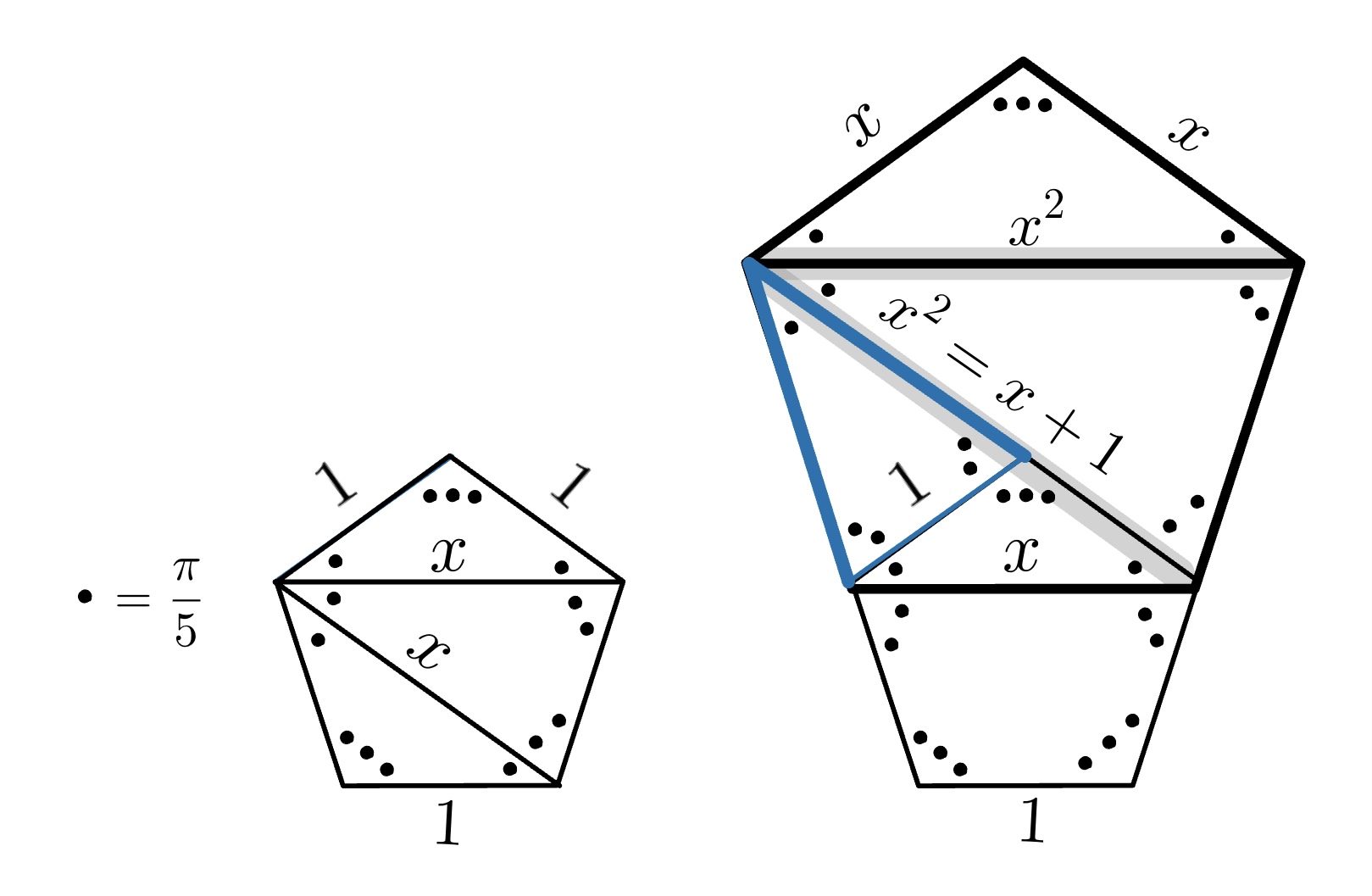} 
    \hbox{\hphantom{dddddddd}}
    \caption{A unit regular pentagon with diagonal length $x$ (left) and two stacked regular pentagons, the upper, larger one of length $x$ and diagonal length $x^2$ (right).}
    \label{figure1}
\end{figure}

We demonstrate the well-known fact
that the diagonal of a unit regular pentagon
has length the golden ratio.
Figure~\ref{figure1} (left) shows a regular unit pentagon, with diagonal length $x$. With hindsight, we know that $x=(\sqrt5 +1)/2$, i.e., the {\it golden ratio.}
However, this may not be known to the reader at first. Each dot represents the same angle. The regularity and consequent symmetry of the pentagon makes clear that three dots can represent the common angle of the pentagon. We see that a single dot thus represents an angle of 36 degrees, or $\pi/5$ radians. 

%\begin{figure}[htbp]
 % \begin{center}
  %  \includegraphics[scale=.5]{PWW.jpg}
    %\includegraphics[trim={20 355 20 160},clip,scale=.75]{IMG_0085.jpg}
  %\end{center}
  %\caption{The diagonal of a regular pentagon.}
  %\label{figure1}
%\end{figure}

The key fact is that $x$ is not only the length of the unit diagonal pentagon, but it is also the ratio of {\it diagonal to length of any regular pentagon}.
Figure~\ref{figure1} (right) has the original unit pentagon, but also with a second regular pentagon of length $x$, built upon the diagonal from the original pentagon. 
Thus, this larger pentagon has diagonal length $x^2$.  Such a diagonal, horizontal and near the top of the figure, is shaded in gray.
Each dot again represents $\pi/5$ radians, and the angles shown are immediate to verify. A second diagonal, also of length $x^2$, is also shown. It has its lower right portion in common with the original unit pentagon (and thus this portion has length $1$). The other portion of this diagonal has length $x$ as it is part of one of the isosceles triangles.

This shows that $x^2=x+1$,
which is famously known to be equivalent to the fact that $x$ is the golden ratio (for example, see~\cite{Livio}, page 80).
The value of this approach is not in the use of isosceles triangles, but in the act of the stacking of the pentagons. This directly introduces $x^2$ into the figure, resulting in a fully visually evident demonstration of all of the relevant parts of the proof.


\bigskip

\begin{thebibliography}{30}

\bibitem{Livio}
  M.~Livio,
  {\em The Golden Ratio: The Story of Phi, the world’s most astonishing number},
  Broadway Books (2002).

\end{thebibliography}
\end{document}